\documentclass{proc-l}

\newtheorem{defn0}{Definition}[section]
\newtheorem{prop0}[defn0]{Proposition}
\newtheorem{conj0}[defn0]{Conjecture}
\newtheorem{thm0}[defn0]{Theorem}
\newtheorem{lem0}[defn0]{Lemma}
\newtheorem{corollary0}[defn0]{Corollary}
\newtheorem{example0}[defn0]{Example}
\newtheorem{remark0}[defn0]{Remark}
\newtheorem{question0}[defn0]{Question}

\newenvironment{defn}{\begin{defn0}}{\end{defn0}}
\newenvironment{prop}{\begin{prop0}}{\end{prop0}}

\newenvironment{thm}{\begin{thm0}}{\end{thm0}}

\newenvironment{cor}{\begin{corollary0}}{\end{corollary0}}
\newenvironment{exm}{\begin{example0}\rm}{\end{example0}}
\newenvironment{rem}{\begin{remark0}\rm}{\end{remark0}}

\newcommand{\defref}[1]{Definition~\ref{#1}}

\newcommand{\thmref}[1]{Theorem~\ref{#1}}

\newcommand{\corref}[1]{Corollary~\ref{#1}}

\newcommand{\secref}[1]{Section~\ref{#1}}
\newcommand{\remref}[1]{Remark~\ref{#1}}

\newcommand{\arrow}[1]{\stackrel{#1}{\longrightarrow}}

\newcommand\Fdot{\mathbf F^\bullet}

\newcommand\Pdot{\mathbf P^\bullet}

\newcommand\Adot{\mathbf A^\bullet}

\newcommand\dm{\operatorname{dim}}

\newcommand\tih{\widetilde H}

\newcommand\mfo{F_{f, \mathbf 0}}

\title{L\^e Modules and Traces}

\author{David B. Massey}

\address{Dept. of Mathematics, Northeastern University, Boston, MA, 02115}

\date{}

\commby{Ronald Fintushel}

\begin{document}





\begin{abstract}  We show how some of our recent results clarify the relationship between the L\^e numbers and the cohomology of the Milnor fiber of a non-isolated hypersurface singularity. The L\^e numbers are actually the ranks of the free Abelian groups -- the L\^e modules -- appearing in a complex whose cohomology is that of the Milnor fiber. Moreover, the Milnor monodromy acts on the L\^e module complex, and we describe the traces of these monodromy actions in terms of the topology of the critical locus.
\end{abstract}

\maketitle

\sloppy

\section{Introduction} \label{sec:intro}

\smallskip

From 1988 to 1995, we developed the L\^e numbers as a generalization to non-isolated affine hypersurface singularities  of the Milnor number of an isolated hypersurface singularity; see \cite{lecycles}. Since producing the L\^e cycles and L\^e numbers, we have wondered how to generalize their construction to the completely general case of a function $f:X\rightarrow\mathbb C$, where $X$ itself may be an arbitrarily singular space and $f$ has an arbitrary critical locus. Of course, what the ``critical locus'' of such an $f$ should mean is part of the problem.

\smallskip

Recently, in \cite{singenrich},  we developed a general framework for producing ``L\^e-like'' cycles, associated to such general $f$. We found that the correct generalization of the L\^e cycles were a type of cycle where the coefficients were not integers, but were, in fact, modules. Moreover, these ``enriched cycles'' were graded. Similarly, the L\^e numbers were replaced by graded modules; to obtain actual numbers, one had to take the Euler characteristics of these graded modules.

The main problem with \cite{singenrich} is that the setting is so general that it is difficult to translate the results into statements that are more down-to-Earth; in particular, it is unclear what the results of \cite{singenrich} tell one about the ``classical'' L\^e numbers of a function on an affine space. We make these translations in \secref{sec:recent} of this paper. 

What we show is that the L\^e numbers are actually the ranks of the free Abelian groups -- the L\^e modules -- appearing in a complex whose cohomology is the integral cohomology of the Milnor fiber. Moreover, the Milnor monodromy acts naturally on this L\^e module complex.

\smallskip

In \secref{sec:traces}, we prove that the traces of the monodromy actions on the L\^e module complex depend solely on the topology/geometry of the critical locus of $f$. This places some interesting restrictions on the characteristic polynomials of the Milnor monodromy action on the cohomology of the Milnor fiber.

\bigskip

\section{Quick Review of L\^e Numbers} \label{sec:review}

\smallskip

Let $\mathcal U$ be an open neighborhood of the origin in $\mathbb C^{n+1}$ and let $f:(\mathcal U, \mathbf 0)\rightarrow(\mathbb C, 0)$ be complex analytic. Let $s:=\dim_{\mathbf 0}\Sigma f$. We will assume throughout this paper that $\mathcal U$ is small enough so that $\Sigma f\subseteq V(f)$.

If $s=0$, then the Milnor number, $\mu_f(\mathbf 0)$, of $f$ at the origin is an extremely important piece of numerical data related to the topology of the hypersurface $V(f)$. In particular, the reduced cohomology, $\tih^*()$, of the Milnor fiber, $F_{f, \mathbf 0}$, of $f$ at the origin is zero in all degrees, except possibly in degree $n$, and $\tih^n(\mfo)\cong\mathbb Z^{\mu_f(\mathbf 0)}$.

\smallskip

Fix a set of linear coordinates $\mathbf z:=(z_0, \dots, z_n)$ for $\mathcal U$. The L\^e numbers of $f$, $\lambda^*_{f, \mathbf z}(\mathbf 0)$, with respect to $\mathbf z$ are a generalization of the Milnor number to the case where $s>0$. See \cite{lecycles}. We shall not recall the definition of the L\^e numbers here, but we will recall some of their properties.

The L\^e numbers are defined algebraically, and their existence requires that the coordinates $\mathbf z$ are chosen generically in a fairly weak sense. If $\mathbf z$ is generic enough so that $\lambda^*_{f, \mathbf z}(\mathbf 0)$ exists, then $\lambda^*_{f, \mathbf z}(\mathbf p)$ exists for all $\mathbf p\in V(f)$ near $\mathbf 0$. It is only the choice of the coordinates $(z_0, \dots, z_{s-1})$ that matter is calculating the L\^e numbers near the origin (recall that $s$ is calculated at the origin).

The $i$-dimensional L\^e number at $\mathbf 0$, $\lambda^i_{f, \mathbf z}(\mathbf 0)$, is possibly non-zero only for those $i$ such that $0\leq i \leq s$. If $s=0$, the only conceivably non-zero L\^e number is $\lambda^0_{f, \mathbf z}(\mathbf 0)$, which equals the Milnor number $\mu_f(\mathbf 0)$. 

If $\mathcal S$ is a stratification of $V(f)$ which contains $V(f)-\Sigma f$, and such that all of the L\^e numbers $\lambda^*_{f, \mathbf z}(\mathbf p)$ remain constant as $\mathbf p$ varies along a stratum, then $\mathcal S$ is a good stratification, i.e., a stratification such that, for all $S\in \mathcal S$, the pair $(\mathcal U-V(f), S)$ satisfies the $a_f$ condition.

If one fixes a good stratification, $\mathcal S$, of $V(f)$ in a neighborhood of the origin,  the coordinates $z_0, \dots, z_{s-1}$ are {\it prepolar} at the origin with respect to $\mathcal S$ provided that, for all $j$ such that $0\leq j\leq s-1$, $V(z_j)$ transversely intersects all of the strata of $\mathcal S\cap V(z_0, \dots, z_{j-1})$ in a neighborhood of the origin except, perhaps, at the origin itself. Prepolar coordinates are generic. 
If $z_0, \dots, z_{s-1}$ are prepolar at the origin with respect to $\mathcal S$, then the coordinates are prepolar with respect to $\mathcal S$ at all points near the origin. If the coordinates $z_0, \dots, z_{s-1}$ are prepolar at the origin with respect to some good stratification, then there is a chain complex 
$$0\rightarrow\mathbb C^{{}^{\lambda^n_f}}\rightarrow\mathbb
C^{{}^{\lambda^{n-1}_f}} \cdots\rightarrow\mathbb C^{{}^{\lambda^{1}_f}}\rightarrow\mathbb C^{{}^{\lambda^0_f}}\rightarrow 0$$  
with
cohomology isomorphic to $\tih^*(F_{f, \mathbf 0};\ \mathbb C)$, where we have written $\lambda^i_f$ for
$\lambda^i_{f, \mathbf z}(\mathbf 0)$, and the $\mathbb C^{{}^{\lambda^i_f}}$ term stands in degree $n-i$ (this complex is mistakenly reversed in \cite{lecycles}). 

\bigskip

\section{More Recent Results} \label{sec:recent}

\smallskip

In \cite{singenrich}, the results that we proved were so general that it may be difficult to see what they say about the L\^e numbers. This section will clarify the matter. Despite the fact that the results here follow quickly from the results of \cite{singenrich}, we state them as theorems, since we have never before presented them as we do below.

\smallskip

We continue to let  $\mathbf z=(z_0, \dots, z_n)$ be a set of linear coordinates for $\mathcal U$, and we will simplify our notation by also writing  $z_0, \dots, z_n$ when we restrict these coordinate functions to various subspaces.

\bigskip

\noindent{\bf Nearby and Vanishing Cycles}

\smallskip

We need to take nearby and vanishing cycles with respect to $f$ and our various coordinates functions. Let us first recall some basic facts about these functors.

Suppose that $h:X\rightarrow\mathbb C$ is an analytic function. Fix a reasonably nice base ring (e.g., a P.I.D. or a field) $R$. Let $D^b_c(X)$ denote the derived category of bounded, constructible complexes of sheaves of $R$-modules on $X$.  

Let $j^*_h$ denote the restriction functor from $D^b_c(X)$ to $D^b_c(V(h))$, which restricts a complex on $X$ to $V(h)$. We wish to discuss two other functors, $\psi_h$ and $\phi_h$, the nearby cycles along $h$ and the vanishing cycles along $h$, respectively, from $D^b_c(X)$ to $D^b_c(V(h))$. 

If $\Adot\in D^b_c(X)$, then the stalk cohomology of $\psi_h\Adot$ at a point $x\in V(h)$ is isomorphic to the hypercohomology of the Milnor fiber of $h$ at $x$ with coefficients in $\Adot$, i.e., $H^i(\psi_h\Adot)_x\cong \mathbb H^i(F_{h, x};\ \Adot)$. In particular, when applied to the constant sheaf, we obtain $H^i(\psi_hR_X^\bullet)_x\cong H^i(F_{h, x};\ R)$.

If $\Adot\in D^b_c(X)$, then the stalk cohomology of $\phi_h\Adot$ at a point $x\in V(h)$ is isomorphic to the shifted relative hypercohomology of the pair consisting of an open ball in $X$ and the Milnor fiber of $h$ at $x$ with coefficients in $\Adot$, i.e., $H^i(\phi_h\Adot)_x\cong \mathbb H^{i+1}({\stackrel{\circ}{B}}_\epsilon(x)\cap X, F_{h, x};\ \Adot)$. In particular, when applied to the constant sheaf, we obtain $H^i(\phi_hR_X^\bullet)_x\cong \widetilde H^i(F_{h, x};\ R)$, where $\widetilde H$ denotes reduced cohomology.

There are natural maps ${\operatorname{comp}}_h: j^*_h\rightarrow \psi_h$ and ${\operatorname{can}}_h: \psi_h\rightarrow \phi_h$; these are the comparison map and canonical map, respectively. For all $\Adot\in D^b_c(X)$, there is the natural distinguished triangle (which we write in-line):
$$
j^*_h\Adot\xrightarrow{{\operatorname{comp}}_h\Adot}\psi_h\Adot\xrightarrow{{\operatorname{can}}_h\Adot}\phi_h\Adot\arrow{[1]}.
$$
The shifted functors $\psi_h[-1]:=\psi_h\circ[-1]$ and $\phi_h[-1]:=\phi_h\circ[-1]$ behave more nicely, in a functorial sense, than do the unshifted functors; in particular, $\psi_h[-1]$ and $\phi_h[-1]$ take perverse sheaves on $X$ to perverse sheaves on $V(h)$. For this reason, we always include the shifts, and we write the canonical distinguished triangle above as
$$
j^*_h[-1]\Adot\xrightarrow{{\operatorname{comp}}_h\Adot}\psi_h[-1]\Adot\xrightarrow{{\operatorname{can}}_h\Adot}\phi_h[-1]\Adot\arrow{[1]},\leqno(\dagger)
$$
where, of course, $j^*_h[-1]:= j^*_h\circ [-1]$.

There are natural Milnor monodromy automorphisms $T_h: \psi_h[-1]\rightarrow \psi_h[-1]$ and $\widetilde T_h: \phi_h[-1]\rightarrow \phi_h[-1]$ such that, for all $\Adot\in D^b_c(X)$, for all $x\in V(h)$, $T_h\Adot$ and $\widetilde T_h\Adot$ induce the Milnor monodromy actions on both $\mathbb H^*(F_{h, x};\ \Adot)$ and  $\mathbb H^{*}({\stackrel{\circ}{B}}_\epsilon(x)\cap X, F_{h, x};\ \Adot)$, respectively. The triple $(\operatorname{id}, T_h\Adot, \widetilde T_h\Adot)$ acts on the triangle $(\dagger)$ commutatively.

\vskip .3in

Now that we have reviewed nearby and vanishing cycles, we will continue to describe some our more recent results. 

\smallskip

Following Definition 5.1 of \cite{singenrich}, but simplifying our notation for our current setting, we make the following definition.

\medskip

\begin{defn}\label{defn:old2.1} For all $j$ such that $0\leq j\leq n$, define $\Phi^j_{f, \mathbf z}$ to be the complex of sheaves of $\mathbb Z$-modules given by iterating the vanishing and nearby cycles
$$
\Phi^j_{f, \mathbf z}:=\phi_{z_j}[-1]\psi_{z_{j-1}}[-1]\dots\psi_{z_0}[-1]\phi_f[-1]\mathbb Z_{\mathcal U}^\bullet[n+1], 
$$
where, when $j=0$, we mean that $\Phi^0_{f, \mathbf z}:=\phi_{z_0}[-1]\phi_f[-1]\mathbb Z_{\mathcal U}^\bullet[n+1]$.

\smallskip

We say that the coordinates $\mathbf z$ are {\it $\phi_f[-1]\mathbb Z_{\mathcal U}^\bullet[n+1]$-isolating at $\mathbf 0$} if and only if, for all $j$ such that $0\leq j\leq s-1$, $\dm_\mathbf 0\operatorname{supp}\big(\Phi^j_{f, \mathbf z}\big)\leq 0$.
\end{defn}

\bigskip

\begin{rem}\label{rem:iteratevan} Iterated vanishing and nearby cycles, such as we use above, seem to have first appeared in the literature in Sabbah's paper \cite{sabbahprox}.
\end{rem}

\bigskip

As $\mathbb Z_{\mathcal U}^\bullet[n+1]$ is a perverse sheaf, and as shifted nearby and vanishing cycles take perverse sheaves to perverse sheaves, we see that each $\Phi^j_{f, \mathbf z}$ is a perverse sheaf. If the coordinates $\mathbf z$ are $\phi_f[-1]\mathbb Z_{\mathcal U}^\bullet[n+1]$-isolating, then each $\Phi^j_{f, \mathbf z}$ is a perverse sheaf with the origin as an isolated point of its support; this implies that the stalk cohomology at the origin of each $\Phi^j_{f, \mathbf z}$ is concentrated in degree $0$. 

\medskip

Hence, we make the following definition:

\medskip

\begin{defn}\label{defn:old2.2} If the coordinates $\mathbf z$ are $\phi_f[-1]\mathbb Z_{\mathcal U}^\bullet[n+1]$-isolating at the origin, then we define the {\it $j$-dimensional L\^e module of $f$ with respect to $\mathbf z$ at the origin} to be
$$
M^j_{f, \mathbf z}:= H^0\big(\Phi^j_{f, \mathbf z}\big)_\mathbf 0.
$$
\end{defn}

\vskip .3in

The next theorem follows from results in \cite{singenrich}. While we cannot reproduce, in a reasonable amount of space, the proofs  from \cite{singenrich} that we need, we can easily describe the differentials $M^{j}_{f, \mathbf z}\arrow{\partial_{j}} M^{j-1}_{f, \mathbf z}$ which appear below.

Let
$$
\Psi^{j-1}_{f, \mathbf z}:=\psi_{z_{j-1}}[-1]\psi_{z_{j-2}}[-1]\dots\psi_{z_0}[-1]\phi_f[-1]\mathbb Z_{\mathcal U}^\bullet[n+1].
$$
Then the canonical morphism in $(\dagger)$ yields a morphism $\Psi^{j-1}_{f, \mathbf z}\xrightarrow{{\operatorname{can}}_{j-1}} \Phi^{j-1}_{f, \mathbf z}$. Applying the functor $\phi_{z_{j}}[-1]$, we obtain the morphism 
$$\Phi^{j}_{f, \mathbf z}= \phi_{z_{j}}[-1]\Psi^{j-1}_{f, \mathbf z}\xrightarrow{\phi_{z_{j}}[-1]({\operatorname{can}}_{j-1})} \phi_{z_{j}}[-1]\Phi^{j-1}_{f, \mathbf z},$$
and a corresponding morphism on stalk cohomology
$$H^0\big(\Phi^{j}_{f, \mathbf z}\big)_{\mathbf 0}\xrightarrow{\partial^\prime_j} H^0\big(\phi_{z_{j}}[-1]\Phi^{j-1}_{f, \mathbf z}\big)_{\mathbf 0}.$$
Now, if the coordinates $\mathbf z$ are $\phi_f[-1]\mathbb Z_{\mathcal U}^\bullet[n+1]$-isolating at the origin, then there exists a neighborhood $\mathcal W$ of the origin such that $\Phi^{j-1}_{f, \mathbf z}$ is zero when restricted to $\mathcal W-\{\mathbf 0\}$ and, hence,  $\mathbf 0$ is not in the support of $\psi_{z_{j}}[-1]\Phi^{j-1}_{f, \mathbf z}$. Therefore, the long exact sequence on stalk cohomology given by the canonical distinguished triangle $(\dagger)$ yields an isomorphism $H^0\big(\phi_{z_{j}}[-1]\Phi^{j-1}_{f, \mathbf z}\big)_{\mathbf 0}\cong H^0\big(\Phi^{j-1}_{f, \mathbf z}\big)_{\mathbf 0}$; by composing this isomorphism with the map $\partial^\prime_j$ above, we obtain a map $M^j_{f, \mathbf z}\arrow{\partial_{j}} M^{j-1}_{f, \mathbf z}$.

\bigskip

\begin{thm}\label{thm:old2.3} Suppose that the coordinates $\mathbf z$ are $\phi_f[-1]\mathbb Z_{\mathcal U}^\bullet[n+1]$-isolating at the origin. Then,

\smallskip

\noindent i) $M^j_{f, \mathbf z}=0$ for all $j>s$;

\smallskip

\noindent ii) if $s:={\operatorname{dim}}_{\mathbf 0}\Sigma f\geq 1$, $M^s_{f, \mathbf z}= \psi_{z_{s-1}}[-1]\psi_{z_{s-2}}[-1]\dots\psi_{z_0}[-1]\phi_f[-1]\mathbb Z_{\mathcal U}^\bullet[n+1]$;

\smallskip

\noindent iii) for all $j$, $M^j_{f, \mathbf z}$ is free Abelian;

\smallskip

\noindent iv)  there is a complex of $\mathbb Z$-modules
$$
0\arrow{\partial_{s+1}} M^s_{f, \mathbf z}\arrow{\partial_{s}} M^{s-1}_{f, \mathbf z}\arrow{\partial_{s-1}} \cdots \arrow{\partial_{2}} M^1_{f, \mathbf z}\arrow{\partial_{1}} M^0_{f, \mathbf z}\arrow{\partial_{0}} 0
$$
such that, for all $j$, the cohomology $\operatorname{ker} \partial_j/\operatorname{im} \partial_{j+1}$ is isomorphic to $\tih^{n-j}(F_{f, \mathbf 0};\ \mathbb Z)$;

\smallskip

\noindent v) the Milnor monodromy automorphism acts on the complex in iv) in a compatible fashion with the Milnor monodromy action on $\tih^*(F_{f, \mathbf 0};\ \mathbb Z)$;

\smallskip

\noindent vi) the induced Milnor monodromy in v) acts on each $M^j_{f, \mathbf z}$  quasi-unipotently, i.e., the complex roots of the characteristic polynomials are all roots of unity.
\end{thm}

\medskip

\begin{proof} We refer to Section 5 of \cite{singenrich} extensively. The complex $\Adot$ of Section 5 of \cite{singenrich} will here be $\phi_f[-1]\mathbb Z_{\mathcal U}^\bullet[n+1]$. The number $d$ of Section 5 of \cite{singenrich} will here be $s$. We will concentrate  our attention at the origin.

\smallskip

Parts i) and ii) follow immediately from Corollary 5.15 of \cite{singenrich}. Part iii)  follows immediately from Theorem 5.23 and Remark 5.24 of \cite{singenrich}. Part iv)  follows immediately from Theorem 5.18  of \cite{singenrich}. Part v) follows from the naturality of the monodromy action on the vanishing cycles.

We need to show Part vi). Let $\widetilde T_f$ denote the monodromy automorphism on $\phi_f[-1]\mathbb Z_{\mathcal U}^\bullet[n+1]$. In a neighborhood of the origin, the Monodromy Theorem \cite{clemens}, \cite{sgavii1}, \cite{landman}, \cite{lemono1}   implies that there exist $a, b\in\mathbb N$ such that $(\operatorname{id}-\widetilde T_f^a)^b = 0$ in the derived category; one applies the functors $\phi_{z_j}[-1]\psi_{z_{j-1}}[-1]\dots\psi_{z_0}[-1]$ to this equality, and the conclusion follows immediately. 
\end{proof}

\bigskip

\begin{defn}\label{defn:old2.4} We refer to the complex appearing in \thmref{thm:old2.3} as the {\it L\^e module complex} of $f$ at $\mathbf 0$ with respect to $\mathbf z$. We refer to the automorphism, $\alpha_j$, on $M^j_{f, \mathbf z}$ which is induced by the Milnor monodromy as the {\it $j$-dimensional L\^e-Milnor monodromy} (or, simply, {\it LM monodromy}).\end{defn}

\bigskip

We must still establish relationships between $\phi_f[-1]\mathbb Z_{\mathcal U}^\bullet[n+1]$-isolating coordinates and prepolar coordinates, and between the L\^e numbers and the ranks of the L\^e modules.

\bigskip

Let $\operatorname{Bl}_{J(f)}\mathcal U\arrow{\pi}\mathcal U$ denote the blow-up of $\mathcal U$ along the Jacobian ideal, $J(f)$, of $f$. This blow-up naturally sits inside $\mathcal U\times \mathbb P^{n}$. Let $E$ denote the exceptional divisor of this blow-up. Below, $\pi_*$ denotes the proper push-forward of cycles.

\bigskip

\begin{thm}\label{thm:old2.5} The following are equivalent:

\smallskip

\noindent a) The coordinates $\mathbf z$ are $\phi_f[-1]\mathbb Z_{\mathcal U}^\bullet[n+1]$-isolating at the origin;

\smallskip

\noindent b) there exists an open neighborhood $\mathcal W$ of the origin inside $\mathcal U$ such that, for all $j$ such that $0\leq j\leq s-1$, $E$ properly intersects $\mathcal W\times(\mathbb P^j\times\{\mathbf 0\})$ inside $\mathcal U\times\mathbb P^n$ (i.e., the intersection is purely $j$-dimensional) and $$\dm_\mathbf 0\Big(V(z_0, \dots, z_{j-1})\cap \pi\big(E\cap (\mathcal W\times(\mathbb P^j\times\{\mathbf 0\})\big)\Big)\leq 0;$$

\smallskip

\noindent c) there exists an open neighborhood $\mathcal W$ of the origin inside $\mathcal U$ such that, for all $j$ such that $0\leq j\leq n$, $E$ properly intersects $\mathcal W\times(\mathbb P^j\times\{\mathbf 0\})$ inside $\mathcal U\times\mathbb P^n$ and $$\dm_\mathbf 0\Big(V(z_0, \dots, z_{j-1})\cap \pi_*\big(E\cdot (\mathcal W\times(\mathbb P^j\times\{\mathbf 0\})\big)\Big)\leq 0.$$

\vskip .3in

\noindent If the equivalent conditions above hold,  then, for all $j$ such that $0\leq j\leq n$,

\smallskip

\noindent i) the germ of the $j$-dimensional L\^e cycle, $\Lambda^j_{f, \mathbf z}$, is defined and is equal to 
$$\pi_*\big(E\cdot (\mathcal W\times(\mathbb P^j\times\{\mathbf 0\})\big);$$

\smallskip

\noindent ii) the $j$-dimensional L\^e number, $\lambda^j_{f, \mathbf z}(\mathbf 0)$ is defined and is equal to the rank of the $j$-dimensional L\^e module $M^j_{f, \mathbf z}$.
\end{thm}

\medskip

\begin{proof} Except for the equivalence of b) to a) and c), this is Theorem 6.7 of \cite{singenrich}, where $\Fdot := \mathbb Z^\bullet_{\mathcal U}[n+1]$ and $\tilde f := f$. The equivalence of b) to a) and c) follows at once from Theorem 5.10 of \cite{singenrich}. \end{proof}

\bigskip

\begin{rem}\label{rem:old2.6} In fact, most of the results of \cite{singenrich} hold with coefficients in extremely general base rings $R$; in particular, they hold when $R$ a finite field such as $\mathbb Z/p\mathbb Z$. 

Thus, the more general form of \thmref{thm:old2.5} tells us that the coordinates $\mathbf z$ are $\phi_f[-1]\mathbb Z_{\mathcal U}^\bullet[n+1]$-isolating at the origin if and only if they are $\phi_f[-1]R_{\mathcal U}^\bullet[n+1]$-isolating at the origin, and when these equivalent conditions hold, there is a complex
$$
0\rightarrow R^{{}^{\lambda^s_{f, \mathbf z}(\mathbf 0)}}\rightarrow R^{{}^{\lambda^{s-1}_{f, \mathbf z}(\mathbf 0)}}\rightarrow\dots\rightarrow R^{{}^{\lambda^0_{f, \mathbf z}(\mathbf 0)}}\rightarrow 0
$$
whose cohomology is isomorphic to $\widetilde H^*(F_{f, \mathbf 0};\ R)$. In other words, the ``L\^e numbers with $R$ coefficients'' remain the same. 

Of course, by the Universal Coefficient Theorem, 
$$
\widetilde H^j(F_{f, \mathbf 0};\ R)\cong \big(\widetilde H^j(F_{f, \mathbf 0};\ \mathbb Z)\otimes R\big)\ \oplus\ \operatorname{Tor}\big(\widetilde H^{j+1}(F_{f, \mathbf 0};\ \mathbb Z), R\big).
$$
Therefore, by varying $R$ through $\mathbb Z$ and through the fields $\mathbb Z/p\mathbb Z$, we find that, not only do the L\^e numbers put restrictions on the Betti numbers of $\widetilde H^*(F_{f, \mathbf 0};\ \mathbb Z)$, they also place upper-bounds on the number of $p$-torsion direct summands of $\widetilde H^*(F_{f, \mathbf 0};\ \mathbb Z)$.\end{rem}

\bigskip

Finally, below, why see why prepolar coordinates are important.

\bigskip

\begin{thm}\label{thm:2.7} Suppose the coordinates $\mathbf z$ are prepolar for $f$ at the origin. Then, the coordinates are $\phi_f[-1]\mathbb Z_{\mathcal U}^\bullet[n+1]$-isolating at the origin,  $f_0:= f_{|_{V(z_0, \dots, z_{s-1})}}$ has an isolated critical point at the origin, and $\widetilde H^{n-s}(F_{f, \mathbf 0};\ \mathbb Z)\cong\operatorname{ker}\partial_s$ injects into $\widetilde H^{n-s}(F_{f_0, \mathbf 0};\ \mathbb Z)$ by a Milnor monodromy-invariant map. \end{thm}

\medskip

\begin{proof} That the coordinates $\mathbf z$ are $\phi_f[-1]\mathbb Z_{\mathcal U}^\bullet[n+1]$-isolating is immediate from the second statement in Theorem 6.5 of \cite{singenrich}. Note, that in our simple setting, the $a_{f, \Fdot}$ partition referred to in 6.5 of \cite{singenrich} is simply a good partition. See Definition 4.7 of \cite{singenrich}, use that $\Fdot := \mathbb Z^\bullet_{\mathcal U}[n+1]$, and note that the only visible stratum is $\mathcal U-V(f)$.

The remaining statements follow from an inductive application of the main result of L\^e in \cite{leattach}.\end{proof}

\bigskip

\begin{rem}\label{rem:slices} In the theorem above, only the statement about propolar coordinates implying that the coordinates are $\phi_f[-1]\mathbb Z_{\mathcal U}^\bullet[n+1]$-isolating is truly new. The comparison of the Milnor fiber of $f$ and of successive hyperplane slices of $f$ is a now-common technique, which originated with work of Hamm, L\^e, and Teissier in 1973; see, for instance, \cite{hammle}, \cite{leattach}, and \cite{teissiercargese}.
\end{rem}

\bigskip

\section{L\^e-Milnor Traces} \label{sec:traces}

\smallskip

In this section, we will show that the traces of the L\^e-Milnor monodromies (\defref{defn:old2.4})  can easily be described in terms of the topology of $\Sigma f$. We continue with $f$, $\mathbf z$, and $s$ as given in the previous sections.

\bigskip

\begin{defn}\label{defn:old3.1} For all $j$ such that $0\leq j\leq s$, we define the {$j$-dimensional (complex) link of $\Sigma f$ at the origin with respect to $\mathbf z$} to be 
$$
\mathbb L^j_{{}_{\Sigma f, \mathbf z}}:=\stackrel{\circ}{B_\epsilon}\cap\Sigma f\cap V(z_0-a_0)\cap V(z_1-a_1)\cap \dots\cap V(z_{s-j-1}-a_{s-j-1}),
$$
where $0\ll |a_0|\ll |a_1|\ll\dots \ll |a_{s-1}|\ll \epsilon\ll 1$. When $j=s$, we mean, of course, that $\mathbb L^s_{{}_{\Sigma f, \mathbf z}}:=\stackrel{\circ}{B_\epsilon}\cap\Sigma f$. It is immediate that
$$H^k(\mathbb L^j_{{}_{\Sigma f, \mathbf z}}; \ \mathbb Z)\cong H^k(\psi_{z_{s-j-1}}\psi_{z_{s-j-2}}\dots\psi_{z_0}\mathbb Z^\bullet_{\Sigma f})_\mathbf 0.$$
\end{defn}

\bigskip

Note that we refer to $\mathbb L^j_{{}_{\Sigma f, \mathbf z}}$ as the ``$j$-dimensional link'', even when the coordinates are such that $\mathbb L^j_{{}_{\Sigma f, \mathbf z}}$ is not $j$-dimensional. However, generically, our terminology makes sense, since we have:

\bigskip

\begin{prop}\label{prop:old3.2} If the coordinates $\mathbf z$ are $\phi_f[-1]\mathbb Z_{\mathcal U}^\bullet[n+1]$-isolating at the origin, then, for all $j$ such that $0\leq j\leq s$, 
$\mathbb L^j_{{}_{\Sigma f, \mathbf z}}$ is $j$-dimensional.
\end{prop}

\medskip

\begin{proof} If $\mathbf z$ are $\phi_f[-1]\mathbb Z_{\mathcal U}^\bullet[n+1]$-isolating, then the L\^e numbers are defined at origin (\thmref{thm:old2.5} ii). As the L\^e numbers are defined, there is an equality of sets $\displaystyle\Sigma f=\bigcup_{i\leq s}\Lambda^i_{f, \mathbf z}$ and, for each $i$, $V(z_0, \dots, z_{i-1})$ properly intersects $\Lambda^i_{f,\mathbf z}$ at the origin. As $\Lambda^s_{f,\mathbf z}$ is non-empty (or, is not zero), the conclusion follows.\end{proof}

\bigskip

Note that it is definitely {\bf not} true that isolating coordinates yield {\bf pure}-dimensional  $\mathbb L^j_{{}_{\Sigma f, \mathbf z}}$.

\bigskip

So that we do not have to break the statement of the following theorem into two cases, we let $\mathbb L^{-1}_{{}_{\Sigma f, \mathbf z}}:=\emptyset$.

\bigskip

\begin{thm}\label{thm:old3.3}Suppose the coordinates $\mathbf z$ are $\phi_f[-1]\mathbb Z_{\mathcal U}^\bullet[n+1]$-isolating at the origin. Then, for all $j$ such that $0\leq j\leq s$, 
the trace of the $j$-dimensional L\^e-Milnor monodromy, $\alpha_j$, is 
$$(-1)^{n-j}\big(\chi(\mathbb L^{s-j-1}_{{}_{\Sigma f, \mathbf z}})-\chi(\mathbb L^{s-j}_{{}_{\Sigma f, \mathbf z}})\big).
$$
\end{thm}

\medskip

\begin{proof} Throughout this proof, if the Milnor monodromy acts on a complex $\Adot$, we let $\mathcal L_\mathbf 0\{\Adot\}$ denote the Lefschetz number of the Milnor monodromy action on the stalk of $\Adot$ at the origin.

There is the fundamental distinguished triangle
$$
\mathbb Z_{V(f)}^\bullet[n]\rightarrow \psi_f[-1]\mathbb Z_{\mathcal U}^\bullet[n+1]\rightarrow\phi_f[-1]\mathbb Z_{\mathcal U}^\bullet[n+1]\arrow{[1]}.
$$
Restricting to $\Sigma f$, we obtain
$$
\mathbb Z_{\Sigma f}^\bullet[n]\rightarrow \big(\psi_f[-1]\mathbb Z_{\mathcal U}^\bullet[n+1]\big)_{|_{\Sigma f}}\rightarrow\big(\phi_f[-1]\mathbb Z_{\mathcal U}^\bullet[n+1]\big)_{|_{\Sigma f}}\arrow{[1]}.
$$
Now, the support of $\phi_f[-1]\mathbb Z_{\mathcal U}^\bullet[n+1]$ is  $\Sigma f$, and  the monodromy action is the identity on $\mathbb Z_{\Sigma f}^\bullet[n]$. In addition, for all $\mathbf p\in\Sigma f$, A'Campo's result \cite{acamp} implies that the Lefschetz number of the monodromy action on $\psi_f[-1]\mathbb Z_{\mathcal U}^\bullet[n+1]$ at $\mathbf p$ is zero. It follows that, for all $j$ such that $0\leq j\leq s-1$, 
$$\mathcal L_\mathbf 0\big\{\psi_{z_{j}}[-1]\dots\psi_{z_0}[-1]\big(\big(\psi_f[-1]\mathbb Z_{\mathcal U}^\bullet[n+1]\big)_{|_{\Sigma f}}\big)\}=0.$$

Therefore, for all $j$ such that $0\leq j\leq s-1$,
$$
\mathcal L_\mathbf 0\big\{\psi_{z_{j}}[-1]\dots\psi_{z_0}[-1]\phi_f[-1]\mathbb Z_{\mathcal U}^\bullet[n+1]\big\}=$$
$$-\chi_\mathbf 0\big(\psi_{z_{j}}[-1]\dots\psi_{z_0}[-1]\mathbb Z_{\Sigma f}^\bullet[n]\big) = (-1)^{n-j}\chi(\mathbb L^{s-j-1}_{{}_{\Sigma f, \mathbf z}}),
$$
where $\chi_{\mathbf 0}$ denotes the Euler characteristic of the stalk at $\mathbf 0$. 

Now, we find that, for all $j$ such that $1\leq j\leq s$,
$$
\mathcal L_\mathbf 0\big\{\phi_{z_{j}}[-1]\psi_{z_{j-1}}[-1]\dots\psi_{z_0}[-1]\phi_f[-1]\mathbb Z_{\mathcal U}^\bullet[n+1]\big\}\ =
$$
$$\ \mathcal L_\mathbf 0\big\{\psi_{z_{j}}[-1]\dots\psi_{z_0}[-1]\phi_f[-1]\mathbb Z_{\mathcal U}^\bullet[n+1]\big\} + $$
$$\mathcal L_\mathbf 0\big\{\psi_{z_{j-1}}[-1]\dots\psi_{z_0}[-1]\phi_f[-1]\mathbb Z_{\mathcal U}^\bullet[n+1]\big\} =
$$
$$
 = (-1)^{n-j}\chi(\mathbb L^{s-j-1}_{{}_{\Sigma f, \mathbf z}}) + (-1)^{n-j-1}\chi(\mathbb L^{s-j}_{{}_{\Sigma f, \mathbf z}}) =$$
 $$ (-1)^{n-j}\big(\chi(\mathbb L^{s-j-1}_{{}_{\Sigma f, \mathbf z}})-\chi(\mathbb L^{s-j}_{{}_{\Sigma f, \mathbf z}})\big).
$$
When $j=0$, we use A'Campo's result again to obtain
$$
\mathcal L_\mathbf 0\big\{\phi_{z_{0}}[-1]\phi_f[-1]\mathbb Z_{\mathcal U}^\bullet[n+1]\big\}=$$
$$\mathcal L_\mathbf 0\big\{\psi_{z_{0}}[-1]\phi_f[-1]\mathbb Z_{\mathcal U}^\bullet[n+1]\big\}\  +\ \mathcal L_\mathbf 0\big\{\phi_f[-1]\mathbb Z_{\mathcal U}^\bullet[n+1]\big\}\ =
$$
$$
(-1)^{n}\chi(\mathbb L^{s-1}_{{}_{\Sigma f, \mathbf z}})+(-1)^{n+1} = (-1)^{n}\big(\chi(\mathbb L^{s-1}_{{}_{\Sigma f, \mathbf z}})-\chi(\mathbb L^{s}_{{}_{\Sigma f, \mathbf z}})\big).
$$
\end{proof}

\bigskip

\begin{rem}\label{rem:tracethm} Note that the alternating sum of the traces telescopes to yield the basic consequence of A'Campo's result from \cite{acamp}: $\mathcal L_{\mathbf 0}\big\{\phi_f[-1]\mathbb Z^\bullet_{\mathcal U}[n+1]\big\} = (-1)^{n+1}$. 

Thus, in a sense, what \thmref{thm:old2.3} and \thmref{thm:old3.3} tell us is that the L\^e module complex refines the monodromy data given by the Monodromy Theorem and A'Campo's Lefschetz number theorem. More precisely,  \thmref{thm:old2.3} and \thmref{thm:old3.3} use the fact that the Monodromy Theorem and A'Campo's theorem hold at points near the origin, and ``compress'' that data into a statement at the origin itself.

Results of this nature -- that use the monodromy data at nearby points to deduce monodromy information at the origin -- have certainly been studied by other authors; see, for instance, the work of Siersma in \cite{siersmavarlad} and of Dimca in Chapter 6 of \cite{dimcasheaves}.

Finally, related to \remref{rem:old2.6}, we should mention that, since A'Campo's result remains true with arbitrary field coefficients, \thmref{thm:old3.3} is also true modulo $p$, when applied to $\alpha^p_j:=\alpha_j\otimes{\operatorname{id}}_{\mathbb Z/p\mathbb Z}$.
\end{rem}

\medskip

\begin{cor}\label{cor:old3.4} Suppose that, at the origin, $\Sigma f$ is smooth and transversely intersected by $V(z_0, \dots, z_j)$ for all $j$ such that $0\leq j\leq s-1$. Suppose the coordinates $\mathbf z$ are $\phi_f[-1]\mathbb Z_{\mathcal U}^\bullet[n+1]$-isolating at the origin. 

Then, the traces of the L\^e-Milnor monodromies are all zero, except the $s$-dimensional one, which is $(-1)^{n-s-1}$. In particular, $\lambda^j_{f, \mathbf z}(\mathbf 0)\neq 1$ for $j\neq s$.
\end{cor}

\medskip

\begin{proof} The first statement follows immediately from \thmref{thm:old3.3}, since all of the links of $\Sigma f$ are contractible. The second statement follows from the fact that an automorphism of $\mathbb Z$ cannot have trace zero.
\end{proof}

\medskip

\begin{cor}\label{cor:old3.5} Suppose the coordinates $\mathbf z$ are $\phi_f[-1]\mathbb Z_{\mathcal U}^\bullet[n+1]$-isolating at the origin. Then, for all $j$ such that $0\leq j\leq s$,  
$$\big|\chi(\mathbb L^{s-j-1}_{{}_{\Sigma f, \mathbf z}})-\chi(\mathbb L^{s-j}_{{}_{\Sigma f, \mathbf z}})\big|\leq \lambda^j_{f, \mathbf z}(\mathbf 0),$$
with equality holding for a given $j$ if and only if the characteristic polynomial of the $j$-dimensional L\^e-Milnor monodromy is either $(t-1)^{\lambda^j_{f, \mathbf z}(\mathbf 0)}$ or $(t+1)^{\lambda^j_{f, \mathbf z}(\mathbf 0)}$.
\end{cor}

\medskip

\begin{proof} This follows at once from \thmref{thm:old3.3} and the fact the complex eigenvalues of the LM monodromies are roots of unity (from \thmref{thm:old2.3} vi).\end{proof}

\medskip

\begin{cor}\label{cor:old3.6} Suppose that $\Sigma f$ is a set-theoretic local complete intersection at the origin. Suppose the coordinates $\mathbf z$ are $\phi_f[-1]\mathbb Z_{\mathcal U}^\bullet[n+1]$-isolating at the origin, and are also $\mathbb Z_{\Sigma f}^\bullet[s]$-isolating at the origin. 

Then, all of the L\^e-Milnor traces are either $0$ or have the same sign as $(-1)^{n-s-1}$.
\end{cor}

\medskip

\begin{proof} The proof of \thmref{thm:old3.3} shows that the $j$-dimensional LM trace is
$$
-\chi_{\mathbf 0}\big\{\phi_{z_{j}}[-1]\psi_{z_{j-1}}[-1]\dots\psi_{z_0}[-1]\mathbb Z_{\Sigma f}^\bullet[n]\big\} =$$
$$(-1)^{n-s-1}\chi_{\mathbf 0}\big\{\phi_{z_{j}}[-1]\psi_{z_{j-1}}[-1]\dots\psi_{z_0}[-1]\mathbb Z_{\Sigma f}^\bullet[s]\big\}.
$$
However, since $\Sigma f$ is a local complete intersection, $\mathbb Z_{\Sigma f}^\bullet[s]$ is a perverse sheaf. Thus, $$\Pdot:=\phi_{z_{j}}[-1]\psi_{z_{j-1}}[-1]\dots\psi_{z_0}[-1]\mathbb Z_{\Sigma f}^\bullet[s]$$ is a perverse sheaf. As the coordinates are $\mathbb Z_{\Sigma f}^\bullet[s]$-isolating, $\Pdot$ is a perverse sheaf which has the origin as an isolated point in its support; hence, the stalk cohomology of $\Pdot$ at the origin is concentrated in degree zero. The conclusion follows.
\end{proof}

\bigskip

\begin{rem}It is conceivable that $\phi_f[-1]\mathbb Z_{\mathcal U}^\bullet[n+1]$-isolating coordinates are automatically $\mathbb Z_{\Sigma f}^\bullet[s]$-isolating, but we do not see a proof of this fact. 
\end{rem}

\bigskip

\section{Examples} \label{sec:examples}

We will conclude by giving two examples.

\bigskip

\begin{exm}\label{exm:line2} Suppose that, at the origin, $f$ has a smooth one-dimensional critical locus, $f_{|_{V(z_0)}}$ has an isolated critical point, and $V(z_0)$ transversally intersects $\Sigma f$.

Then, $\lambda^1_{f, \mathbf z}(\mathbf 0)$ is simply equal to the Milnor number, $\stackrel\circ\mu$, of a generic hyperplane slice of $f$ at a point on $\Sigma f$, sufficiently close to the origin. The $0$-dimensional L\^e number, $\lambda^0_{f, \mathbf z}(\mathbf 0)$, is defined to be the intersection number $\left(\Gamma^1_{f, z_0}\cdot V\left(\frac{\partial f}{\partial z_0}\right)\right)_{\mathbf 0}$, where $\Gamma^1_{f, z_0}$ is the relative polar curve of $f$ with respect to $z_0$.

\bigskip

Up to isomorphism, the L\^e module complex of \thmref{thm:old2.3} becomes
$$
0\rightarrow {\mathbb Z}^{\lambda^1_{f, \mathbf z}(\mathbf 0)}\arrow{\partial_1} {\mathbb Z}^{\lambda^0_{f, \mathbf z}(\mathbf 0)}\rightarrow 0,
$$
where $\operatorname{ker}{\partial_1}\cong\widetilde H^{n-1}(F_{f, \mathbf 0})$ and $\operatorname{coker}{\partial_1}\cong\widetilde H^{n}(F_{f, \mathbf 0})$. \thmref{thm:old3.3} tells us that the $0$-dimensional LM monodromy has trace $0$, while the $1$-dimensional LM monodromy has trace $(-1)^{n-2}=(-1)^n$.

In addition, by \remref{rem:old2.6} and \remref{rem:tracethm}, the modulo $p$ version of the previous paragraph holds: for each prime $p$, we have a complex 
$$
0\rightarrow \left({\mathbb Z/p\mathbb Z}\right)^{\lambda^1_{f, \mathbf z}(\mathbf 0)}\arrow{\partial^p_1} \left({\mathbb Z/p\mathbb Z}\right)^{\lambda^0_{f, \mathbf z}(\mathbf 0)}\rightarrow 0,
$$
where $\operatorname{ker}{\partial^p_1}\cong\widetilde H^{n-1}(F_{f, \mathbf 0};\ \mathbb Z/p\mathbb Z)$ and $\operatorname{coker}{\partial^p_1}\cong\widetilde H^{n}(F_{f, \mathbf 0};\ \mathbb Z/p\mathbb Z)$. The modulo $p$  LM-monodromies, $\alpha^p_0$ and $\alpha^p_1$, have traces equal to $0$ and  $(-1)^n$ in $\mathbb Z/p\mathbb Z$, respectively.

\bigskip

 If $\lambda^1_{f, \mathbf z}(\mathbf 0) =1$, then a well-known result of Siersma from \cite{siersmaisoline} tells us that either $\lambda^0_{f, \mathbf z}(\mathbf 0)=0$, or $\partial_1\neq 0$ and, for all primes $p$, $\partial^p_1\neq 0$. It follows that there are two cases. In the first case, $\lambda^0_{f, \mathbf z}(\mathbf 0)=0$, $\widetilde H^{n-1}(F_{f, \mathbf 0})\cong {\mathbb Z}$, and $\widetilde H^{n}(F_{f, \mathbf 0})=0$. In the second case, $\lambda^0_{f, \mathbf z}(\mathbf 0)\neq 0$, $\widetilde H^{n-1}(F_{f, \mathbf 0})=0$, and $\widetilde H^{n}(F_{f, \mathbf 0})\cong {\mathbb Z}^{\lambda^0_{f, \mathbf z}(\mathbf 0)-1}$

\bigskip

What can one conclude if $\lambda^1_{f, \mathbf z}(\mathbf 0)=2$? A recent result of ours with L\^e \cite{lemassey} tells us that, when $\dm_{\mathbf 0}\Sigma f =1$, then the above statement -- that either $\lambda^0_{f, \mathbf z}(\mathbf 0)=0$, or $\partial_1\neq 0$ and, for all primes $p$, $\partial^p_1\neq 0$ -- is always true.

As the LM-monodromies act invariantly on the kernel and image of $\partial_1$ (resp., $\partial_1^p$), and as $\operatorname{trace}(\alpha_1)=(-1)^n$ (resp., $\operatorname{trace}(\alpha^p_1)=(-1)^n$ in $\mathbb Z/p\mathbb Z$), we cannot have that the rank (resp., dimension) of the kernel of $\partial_1$ (resp., $\partial_1^p$) is equal to $1$.

Thus, it follows once again that there are two cases. In the first case, $\lambda^0_{f, \mathbf z}(\mathbf 0)=0$, $\widetilde H^{n-1}(F_{f, \mathbf 0})\cong {\mathbb Z}^2$, and $\widetilde H^{n}(F_{f, \mathbf 0})=0$. In the second case, $\lambda^0_{f, \mathbf z}(\mathbf 0)\geq 2$, $\widetilde H^{n-1}(F_{f, \mathbf 0})=0$, and $\widetilde H^{n}(F_{f, \mathbf 0})\cong {\mathbb Z}^{\lambda^0_{f, \mathbf z}(\mathbf 0)-2}$

\end{exm}

\bigskip

\begin{exm}\label{exm:cone} Let us look at an example in which the critical locus has more interesting topology. Suppose that, at the origin, $f$ has a critical locus which is  the $2$-dimensional cone $V(x^2+y^2+z^2, w_1, \dots, w_{n-2})$, where our coordinates on affine space are $\mathbf z:=(x, y, z, w_1, \dots, w_{n-2})$. We shall see that the topology of $\Sigma f$ forces the lower-dimensional L\^e numbers to be non-zero, and puts other constraints on the cohomology and monodromy of the Milnor fiber of $f$ at the origin.

We denote homotopy-equivalence by $\simeq$. Of course,  ${\mathbb L}^2_{{}_{\Sigma f, \mathbf z}}\simeq point$. The $1$-dimensional link, ${\mathbb L}^1_{{}_{\Sigma f, \mathbf z}}$, is isomorphic to the Milnor fiber of $y^2 +z^2$in $\mathbb C^2$ and, hence, ${\mathbb L}^1_{{}_{\Sigma f, \mathbf z}}\simeq S^1$. One easily finds that ${\mathbb L}^0_{{}_{\Sigma f, \mathbf z}}\simeq 2 \ points$. Therefore, the  respective Euler characteristics of  ${\mathbb L}^2_{{}_{\Sigma f, \mathbf z}}$, ${\mathbb L}^1_{{}_{\Sigma f, \mathbf z}}$, and  ${\mathbb L}^0_{{}_{\Sigma f, \mathbf z}}$ are $1$, $0$, and $2$.

Assuming that the coordinates are  $\phi_f[-1]\mathbb Z_{\mathcal U}^\bullet[n+1]$-isolating at the origin, we find that $\operatorname{trace}(\alpha_2) = (-1)^{n-2}(0-2) = (-1)^{n-1}2$, $\operatorname{trace}(\alpha_1) = (-1)^{n-1}(2-0) = (-1)^{n-1}2$, and $\operatorname{trace}(\alpha_0) = (-1)^{n-2}(0-1) = (-1)^{n-1}$.  Furthermore, Corollary 3.3 of \cite{lemassey} tells us that $\partial_2$  cannot be the zero map. 

Thus,  we conclude from \corref{cor:old3.5} that $\lambda^2_{f, \mathbf z}(\mathbf 0)$ and $\lambda^1_{f, \mathbf z}(\mathbf 0)$ much each be at least $2$, and $\lambda^0_{f, \mathbf z}(\mathbf 0)\geq 1$. We also know that $\widetilde H^{n-2}(F_{f, \mathbf 0})$ is zero or isomorphic to $\mathbb Z$.

\bigskip

For special values of the L\^e numbers, one can reach further conclusions. As in the previous example, let $\stackrel\circ\mu$ denote the Milnor number of a normal slice (i.e., the restriction to a generic codimension $2$ linear subspace) of $f$ at a generic point on $\Sigma f$, sufficiently close to the origin. Then, $\lambda^2_{f, \mathbf z}(\mathbf 0)=\stackrel\circ\mu\cdot\left(\Sigma f\cdot V(x, y)\right)_{\mathbf 0} = 2\stackrel\circ\mu$, and so we see that $\lambda^2_{f, \mathbf z}(\mathbf 0)$ obtains its minimum value of $2$ precisely when the Milnor number of a normal slice at a nearby generic point on $\Sigma f$ is $1$. This is the case, for instance, for $f_1:= w_1^2 +(x^2+y^2+z^2)^2$.

\medskip

Let us assume now that $\stackrel\circ\mu=1$ and, hence, $\lambda^2_{f, \mathbf z}(\mathbf 0)=2$.  What more can we say about this case? 

\medskip

Since $\operatorname{trace}(\alpha_2) = (-1)^{n-1}2$,  \corref{cor:old3.5} tells us that the characteristic polynomial of $\alpha_2$, ${\operatorname{char}}_{\alpha_2}(t)$,  would have to equal $(t-(-1)^{n-1})^2$. If $\widetilde H^{n-2}(F_{f, \mathbf 0})=0$, then $(t-(-1)^{n-1})^2$ must divide ${\operatorname{char}}_{\alpha_1}(t)$. If $\widetilde H^{n-2}(F_{f, \mathbf 0})\cong\mathbb Z$, then at least $(t-(-1)^{n-1})$ must divide ${\operatorname{char}}_{\alpha_1}(t)$.

Now, if  $\lambda^1_{f, \mathbf z}(\mathbf 0)=2$, then ${\operatorname{char}}_{\alpha_1}(t)$ must also be $(t-(-1)^{n-1})^2$. If $\lambda^1_{f, \mathbf z}(\mathbf 0)=3$, we cannot be in the case where $\widetilde H^{n-2}(F_{f, \mathbf 0})=0$, i.e., the case where $\operatorname{im}\partial_2 \cong \mathbb Z^2$; for, otherwise, $\operatorname{trace}(\alpha_2)  = (-1)^{n-1}2$ and $\operatorname{trace}(\alpha_1)  = (-1)^{n-1}2$ would imply that $0$ would have to be the ``third'' eigenvalue of the isomorphism $\alpha_1$. Thus, if $\lambda^1_{f, \mathbf z}(\mathbf 0)=3$, then $\operatorname{im}\partial_2 \cong \mathbb Z$, and the restriction of $\alpha_1$ to $\operatorname{im}\partial_2$ has characteristic polynomial $(t-(-1)^{n-1})$; it follows that ${\operatorname{char}}_{\alpha_1}(t)= (t-(-1)^{n-1})(t^2-(-1)^{n-1}t+1)$, where it is important to note that $(t^2-(-1)^{n-1}t+1)$ is irreducible over $\mathbb Z$ (or $\mathbb Q$).

Therefore, when $\Sigma f$ is a $2$-dimensional cone and $\stackrel\circ\mu=1$, we have the following cases:

\vskip .1in

\noindent $a)$\hskip .1in $\lambda^1_{f, \mathbf z}(\mathbf 0)=2$, $\widetilde H^{n-2}(F_{f, \mathbf 0})=0$, $\operatorname{rank}\widetilde H^{n-1}(F_{f, \mathbf 0})=0$, and $\widetilde H^{n}(F_{f, \mathbf 0})\cong\mathbb Z^{\lambda^0_{f, \mathbf z}(\mathbf 0)}$.

\bigskip

\noindent $b)$\hskip .1in $\lambda^1_{f, \mathbf z}(\mathbf 0)=2$, $\widetilde H^{n-2}(F_{f, \mathbf 0})\cong \mathbb Z$, $\operatorname{rank}\widetilde H^{n-1}(F_{f, \mathbf 0})=0$, and $\operatorname{rank}\widetilde H^{n}(F_{f, \mathbf 0})=\lambda^0_{f, \mathbf z}(\mathbf 0)-1$.

\bigskip

\noindent $c)$\hskip .1in $\lambda^1_{f, \mathbf z}(\mathbf 0)=2$, $\widetilde H^{n-2}(F_{f, \mathbf 0})\cong \mathbb Z$, $\operatorname{rank}\widetilde H^{n-1}(F_{f, \mathbf 0})=1$, and $\widetilde H^{n}(F_{f, \mathbf 0})\cong{\mathbb Z}^{\lambda^0_{f, \mathbf z}(\mathbf 0)}$.

\bigskip

\noindent $d)$\hskip .1in $\lambda^1_{f, \mathbf z}(\mathbf 0)=3$, $\widetilde H^{n-2}(F_{f, \mathbf 0})\cong \mathbb Z$, $\operatorname{rank}\widetilde H^{n-1}(F_{f, \mathbf 0})=0$, and $\operatorname{rank}\widetilde H^{n}(F_{f, \mathbf 0})=\lambda^0_{f, \mathbf z}(\mathbf 0)-2$.

\bigskip

\noindent $e)$\hskip .1in $\lambda^1_{f, \mathbf z}(\mathbf 0)=3$, $\widetilde H^{n-2}(F_{f, \mathbf 0})\cong \mathbb Z$, $\operatorname{rank}\widetilde H^{n-1}(F_{f, \mathbf 0})=2$, the monodromy on the free part of $\widetilde H^{n-1}(F_{f, \mathbf 0})$ has characteristic polynomial $(t^2-(-1)^{n-1}t+1)$, and $\widetilde H^{n}(F_{f, \mathbf 0})\cong{\mathbb Z}^{\lambda^0_{f, \mathbf z}(\mathbf 0)}$.

\bigskip

\noindent $f)$\hskip .1in $\lambda^1_{f, \mathbf z}(\mathbf 0)\geq 4$.

\bigskip

As one should expect, we see that knowing the traces of the LM-monodromies, and knowing that the eigenvalues are roots of unity, yield strong results for small values of the L\^e numbers.

\end{exm}

\bibliographystyle{plain}
\bibliography{Masseybib}
\end{document}